\documentclass{amsproc}

\usepackage[margin=1.2in,nomarginpar]{geometry}
\usepackage{amssymb}
\usepackage{amsmath}
\usepackage{mathdots}
\usepackage{amsbsy}
\usepackage{amscd}
\usepackage{amsthm}

\usepackage{url}
\usepackage{xcolor}
\usepackage{amsmath,amssymb,amsthm,amsfonts,longtable}
\usepackage[english]{babel}
\usepackage{tikz-cd}
\usetikzlibrary{cd}
\usetikzlibrary{decorations.markings}
\tikzset{negated/.style={
		decoration={markings,
			mark= at position 0.5 with {
				\node[transform shape] (tempnode) {$\times$};
			}
		},
		postaction={decorate}
	}
}
\usepackage{array}
\usepackage[colorlinks,citecolor=red,urlcolor=blue,bookmarks=false,hypertexnames=true]{hyperref}

\newtheorem{theorem}{Theorem}
\newtheorem{corollary}[theorem]{Corollary}
\newtheorem{proposition}[theorem]{Proposition}
\newtheorem{lemma}[theorem]{Lemma}
\newtheorem{remark}[theorem]{Remark}
\newtheorem{consequence}[theorem]{Consequence}
\newcommand{\restr}{\mathord\downarrow} 
\newcommand{\ind}{\mathord\uparrow} 
\newtheorem{example}[theorem]{Example}
\newtheorem{definition}{Definition}[subsection]
\newcommand{\Irr}{\textnormal{Irr}}
\newcommand{\cd}{\textnormal{cd}}
\newcommand{\nl}{\textnormal{nl}}
\newcommand{\lin}{\textnormal{lin}}

\newcommand{\Core}{\textnormal{Core}}
\newcommand{\inertiagroup}{\textnormal{I}}

\title[Faithful quasi-permutation representations]{Minimal Faithful Quasi-Permutation Representation Degree\\ of $p$-Groups with Cyclic Center}
\author{Sunil Kumar Prajapati$^*$}
\address{Indian Institute of Technology, Bhubaneswar, Arugul Campus, Jatni, Khurda-752050, India.}
\email{skprajapati@iitbbs.ac.in}
\author{Ayush Udeep}
\address{Indian Institute of Technology, Bhubaneswar, Arugul Campus, Jatni, Khurda-752050, India.}
\email{udeepayush@gmail.com}
\thanks{$^{\textbf{*}}$Corresponding author.
}
\subjclass[2010]{primary 20D15; secondary 20C15, 20B05}
\keywords{$p$-groups, Metabelian $p$-groups, Cyclic center, Quasi-permutation representations}
\begin{document}
	\maketitle
	
	\begin{abstract}
		For a finite group $G$, we denote by $\mu(G)$, and $c(G)$, the minimal degree of faithful permutation representation of $G$, and the minimal degree of faithful representation of $G$ by quasi-permutation matrices over the complex field $\mathbb{C}$, respectively. In this article, we study $\mu(G)$, and $c(G)$ for various classes of finite non-abelian $p$-groups with cyclic center. We prove a result for normally monomial $p$-groups with cyclic center which generalizes a result of Behravesh for finite $p$-groups of nilpotency class 2 with cyclic center \cite[Theorem 4.12]{HB}. We also compute minimal degrees for some classes of metabelian $p$-groups.
	\end{abstract}
	
	\section{Introduction}
	Throughout this paper, $G$ is a finite group and $p$ is a prime. Cayley's theorem states that $G$ can be embedded into $S_{|G|}$, the symmetric group on $|G|$ symbols. The minimal faithful permutation degree $\mu(G)$ of $G$ is the least positive integer $n$ such that $G$ is isomorphic to some subgroup of $S_{n}$.
	Wong \cite{W} defined a quasi-permutation group, similar to the definition of a permutation group, as a finite linear group $G$ of degree $n$ such that the trace of every element of $G$ is a non-negative integer.
	This leads to another degree $c(G)$, which denotes the minimal degree of a faithful representation of $G$ by complex quasi-permutation matrices (square matrices over complex field with non-negative integral trace) (see \cite{BGHS}).
	Since every permutation matrix is a quasi-permutation matrix, it is easy to see that $c(G)\leq \mu(G)$, i.e., $c(G)$ provides a lower bound for $\mu(G)$. 
	In recent years, several researchers studied $\mu(G)$ and $c(G)$ extensively (see \cite{AB,  HB, HB1997, BG, BGHS, GA, DLJ, SA, DW}). In \cite[Theorem 3.2]{BG}, Behravesh and Ghaffarzadeh proved that if $G$ is a finite $p$-group of odd order, then $c(G)=\mu(G)$. In \cite{HB}, Behravesh proved the following result.
	
	\begin{theorem} \cite[Theorem 4.12]{HB}\label{thm_HB_class2} If $G$ is a $p$-group of nilpotency class 2 with cyclic center, then $c(G) = |G/Z(G)|^{1/2} |Z(G)|$. 
	\end{theorem} 
	\noindent In \cite{HBmetacyclic}, Behravesh computed $c(G)$, where $G$ is a metacyclic $2$-group with cyclic center. In literature, there are no results available about $c(G)$ for a general $p$-group with cyclic center.
	In this article, we study $\mu(G)$ and $c(G)$ for various classes of finite $p$-groups of odd order with cyclic center. We begin with finding a range for $c(G)$ for various classes of finite non-abelian $p$-groups.\\
	
	\noindent  In \cite[Lemma 2.2]{BG}, Behravesh and Ghaffarzadeh gave an algorithm for the computation of $c(G)$. Let $X \subset 
	\Irr(G)$ such that $\cap_{\chi \in X} \ker (\chi)= 1$ 
	and $\cap_{\chi \in Y} \ker (\chi) \neq 1$ 
	for every proper subset $Y$ of $X$. 
	Let $\xi_X = \sum_{\chi \in X} \left[ \sum_{\sigma \in \Gamma(\chi)}
	\chi^{\sigma}  \right]$ and let $m(\xi_X)$ be the absolute value of 
	the minimum value that $\xi_X$ takes over $G$.
	Then $$c(G) = \min \{\xi_X(1) + m(\xi_X) \; | \; X \subset \Irr(G) 
	\text{\ satisfying\ the\ above\ property} \}.$$
	We call a set $X_G\subset \Irr(G)$, a minimal faithful quasi-permutation representation of $G$ if
	\begin{equation}\label{eq:X_G} 
		\bigcap_{\chi \in X_G} \ker (\chi) = 1 \text{ and } 	\bigcap_{\chi \in Y} \ker (\chi) \neq 1 \text{ for any } Y\subset X_G,
	\end{equation}
	such that $c(G) = \xi(1) + m(\xi)$. In general, $X_{G}$ may contain linear characters of $G$. In Theorem \ref{lemmalinearcharacter}, we give a class of non-abelian $p$-groups in which the linear characters do not belong in $X_{G}$.
	\begin{theorem}	\label{lemmalinearcharacter}
		Let $G$ be a non-abelian $p$-group with $d(Z(G)  \cap G') = d(Z(G))$. 	Then $X_{G} \cap \lin(G) = \emptyset$ and $p^{s+1}$ divides $c(G)$, where $p^s = \min\{ \chi(1)~|~ \chi \in \nl(G) \}$.
	\end{theorem}
		
	\noindent  Now if $G$ is a finite non-abelian $p$-group, then we can see that $c(G) = f(p) = a_1 p + a_2 p^2 + \cdots + a_{r}p^{r}$, and $0 \leq a_{i} < p$ such that $\sum_{i=1}^{r}a_{k} = d(Z(G))$. In the view of the above fact, we prove the following result for any finite non-abelian $p$-group. 
	
	\begin{theorem} \label{P2}
		Let $G$ be a non-abelian $p$-group such that $\exp(G) = p^{b}$, and $\max\cd(G) = p^{e}$. Suppose $c(G) = f(p) = a_1 p + a_2 p^2 + \cdots + a_{r}p^{r}$, and $0 \leq a_{i} < p$. Then $b \leq r$. Further, if $d(Z(G))<p$, then $b \leq r \leq b+e$. The bounds are best possible.
	\end{theorem}
	\noindent From  Theorem \ref{P2}, we deduce the following  corollary. 
	\begin{corollary} 
		\label{cor:newrange2}
		Let $G$ be a non-abelian $p$-group such that $\exp(G) = p$, and $\cd(G) = \{ 1, p \}$. Then $c(G) \leq d(Z(G))p^{2}$. Further, if $p\geq 3$ and $G$ is not a direct product of an abelian and a non-abelian group then $c(G) = d(Z(G))p^2$.
	\end{corollary}
	
	\noindent Theorem \ref{thm:cyclicexistence} guarantees the existence of a non-abelian permutation group of order $p^n$ ($p\geq 3$) acting transitively on $p^{i}$ points, for each $i$ ($2 \leq i \leq n-1$), where $p\geq n$.
	\begin{theorem} 
		\label{thm:cyclicexistence}
		Let $G$ be a non-abelian $p$-group of order $p^{n} ~(p\geq 3)$ with cyclic center. Then 
		\[ c(G) = \mu(G) \in \{ p^2, p^3, \ldots, p^{n-1} \}. \] Moreover, for every $n\geq 3$, and every $i$ in $\{2,\ldots, n-1\}$, there exists a $p$-group of order $p^{n}$, for some prime $p$, say $G_{i}$, such that $\mu(G_{i}) = p^{i}$.
	\end{theorem}
	\noindent  Suppose  $G$ is a non-abelian $p$-group of order $p^{n}$ ($p\geq 3$) with cyclic center. Then from Theorem \ref{thm:cyclicexistence}, $p^2$ divides $c(G)$ and $c(G)$ divides $p^{n-1}$. In Corollary \ref{cor:newrangecyclic}, we improve these bounds with the help of Theorem \ref{P2}.
	
	\begin{corollary} 
		\label{cor:newrangecyclic}
		Let $G$ be a non-abelian $p$-group with cyclic center. Suppose $p^{e} = \max \cd(G)$ and $p^{\alpha} = \min\{ \chi(1)~|~ \ker(\chi) = 1, \chi\in \nl(G) \}$. Then
		\[ p^{\alpha}|Z(G)| \text{ divides } c(G) \text{ and } c(G) \text{ divides } p^{e}\cdot \exp(G). \]
		In particular, if $\cd(G) = \{ 1, p, p^{e} \} ~(e>1)$ and $\exp(G) = p$, then $c(G) = p^{e+1}$.
	\end{corollary}
	
	\noindent  A group $G$ is called a normally monomial group if every complex irreducible character of $G$ is induced from a linear character of a normal subgroup of $G$. Note that the $p$-groups satisfying the hypothesis of Corollary \ref{cor:newrangecyclic} may not be normally monomial. In Theorem \ref{lemma:normallymonomial}, we improve the bounds of $c(G)$ obtained in Corollary \ref{cor:newrangecyclic} for normally monomial $p$-groups. 
	
	\begin{theorem}	\label{lemma:normallymonomial}
		Let $G$ be a normally monomial $p$-group with cyclic center. Suppose $A$ is an abelian normal subgroup of maximum order in $G$. Then 
		\[ (\max \cd(G))|Z(G)| \text{ divides } c(G) \text{ and } c(G)\text{ divides } (\max \cd(G))\exp(A). \]
	\end{theorem}
	
	\begin{remark}
		\textnormal{Let $G$ be a normally monomial $p$-group with cyclic center. Suppose $A$ is an abelian normal subgroup of maximum order in $G$ such that $|Z(G)| = \exp(A)$. Then from Theorem \ref{lemma:normallymonomial}, we get $c(G) = (\max \cd(G))|Z(G)|$. If $G$ is a group of nilpotency class 2, then $G$ is a metabelian group, and hence, a normally monomial group. Further, in Proposition \ref{lemma:class2}, we prove that if $G$ is a finite $p$-group ($p\geq 3$) of nilpotency class 2 with cyclic center, then there exists an abelian normal subgroup $A$ of maximum order such that $|Z(G)| = \exp(A)$. Hence, we deduce Theorem \ref{thm_HB_class2} through Theorem \ref{lemma:normallymonomial} for $p\geq 3$.}
	\end{remark}

	\noindent We conclude the article with a result on metabelian $p$-groups. The class of metabelian $p$-groups is large and since computation of $c(G)$ and $\mu(G)$ heavily depends on the structure of the group $G$, it is difficult to compute $c(G)$ and $\mu(G)$ for any general metabelian $p$-group. In Theorem \ref{prop:elementaryabelian}, we step in this direction by dealing with some classes of metabelian $p$-groups. 
	\begin{remark} \label{remark:nMcd(G)1,p}
		\textnormal{Suppose $G$ is a non-abelian $p$-group ($p\geq 3$) such that $\cd(G) = \{ 1,p \}$, $\exp(G) \in \{ p, p^2 \}$ and $d(Z(G)) < p$. Then from Theorem \ref{P2}, we get $c(G) = ap + bp^2 + cp^3$ for some $0 \leq a, b, c < p$ such that $a+b+c = d(Z(G))$. Further, if $G$ is not a direct product of an abelian and a non-abelian subgroup, then from Lemma \ref{lemma:a_1=0} and \cite[Theorem 3.2]{BG}, we get $c(G) = ap^2 + bp^3$, for some $0 \leq a, b < p$ such that $a+b = d(Z(G))$. We illustrate in Example \ref{example:elementaryabelian}(1) that, this bound is best possible.} 
	\end{remark}
	\noindent In Theorem \ref{prop:elementaryabelian}, under certain conditions we deal with a class of non-abelian $p$-groups such that $\cd(G) = \{ 1,p, p^2 \}$ and find a range for $c(G)$. In fact, we prove the following.

	\begin{theorem}\label{prop:elementaryabelian}
		Let $G$ be a non-abelian $p$-group $(p\geq 3)$ such that  $\exp(G) \in \{ p, p^2 \}$, $\cd(G) = \{ 1, p, p^2 \}$, $d(Z(G)) \geq 2$, and $G$ is not a direct product of an abelian and a non-abelian subgroup. Suppose there exists an elementary abelian normal subgroup of index $p^2$ in $G$.
		Then $c(G) = ap^2 + bp^3$, for some $0 \leq a, b < p$ such that $a+b = d(Z(G))$.
	\end{theorem}

	\noindent	In Section \ref{section:notation}, we summarize the notations and some preliminary results, which are used in the rest of this article. We prove our results in Section \ref{sec:results}. In Subection \ref{subsec:rangeforc(G)}, we compute a range of $c(G)$ for various classes of $p$-groups $G$. 
	In Subsection \ref{subsec:cycliccenter}, we study $c(G)$ for normally monomial $p$-groups with cyclic center. In the same subsection, we compute $\mu(G)$ and $c(G)$, when $G$ is a $p$-group ($p\geq 3$) with cyclic center. In the final subsection, i.e., Subsection \ref{subsec:metabelian}, we present the proof of Theorem \ref{prop:elementaryabelian}.
	
	\section{Notations and Preliminaries}\label{section:notation}
	We denote the empty set by $\emptyset$. For a finite group $G$, 
	$d(G)$ and $G'$ denote the minimal number of generators and the commutator subgroup of $G$, respectively. 
	For $g\in G$, $o(g)$ denotes the order of the element $g$ and for a subgroup $H$ of $G$,
	$\Core_{G}(H)$ denotes the core of $H$ in $G$.
	Let $\Irr(G)$, $\lin(G)$ and $\nl(G)$ be the set of irreducible complex characters, the set of linear characters and the set of nonlinear irreducible characters of $G$, respectively.
	We denote the character degree set, i.e., $\{ \chi(1) ~|~ \chi \in \Irr(G) \}$ of $G$ by $\cd(G)$. For $\chi\in \Irr(G)$, $\mathbb{Q}(\chi)$ denotes the field obtained by adjoining  the values $\chi(g)$, for all
	$g \in G$, to $\mathbb{Q}$, and $\Gamma(\chi)$ denotes the Galois group of $\mathbb{Q}(\chi)$ over $\mathbb{Q}$.
	We denote the elementary abelian $p$-group of rank $k$ by $C_{p}^{k}$, whenever $k\geq 3$.  
	Let	$\phi(n)$ be the Euler phi function, and $\omega_{n}$ be a primitive $n^{th}$ root of unity. \\

 Now, we quote some useful results, which we will use throughout our article. Readers can see \cite{DLJ,DW} for some interesting background results on $\mu(G)$. Johnson \cite{DLJ} has proved the following result, which gives the number of transitive constituents of a minimal faithful permutation representation of a $p$-group. 
	\begin{lemma} \textnormal{\cite[Theorem 3]{DLJ}} \label{L2} 
		Let $G$ be a $p$-group, and let $\mathcal{H}= \{ G_{1},\ldots,G_{n} \}$ be a minimal faithful permutation representation of $G$. If $p$ is odd, then  $n=d(Z(G))$. If $p=2$, then $d(Z(G))/2\leq n\leq d(Z(G))$, the bound $n=d(Z(G))$ being achieved.
	\end{lemma}

	 Let $\chi, ~\psi \in \Irr(G)$. 
	We say that $\chi$ and $\psi$ are Galois conjugate over $\mathbb{Q}$ if there exists $\sigma \in \Gamma(\chi)$ such that $\chi^{\sigma}= \psi$. One can check that Galois conjugacy defines an equivalence relation on $\Irr(G)$. Moreover, if $\mathcal{C}$ denotes the equivalence class of $\chi$ with respect to Galois conjugacy over $\mathbb{Q}$, then $|\mathcal{C}|=| \mathbb{Q}(\chi) : \mathbb{Q} |$ (see \cite[Lemma 9.17]{I}). 
	 For $\psi \in \Irr(G)$, let $\Psi$ be its Galois-sum over $\mathbb{Q}$, i.e., $\Psi = \sum_{\sigma\in \Gamma(\psi)} \psi^{\sigma}$. Let $\{ 1_{G} = \Psi_{0}, \Psi_{1}, \ldots, \Psi_{r} \}$ be all the Galois-sums of irreducible characters of $G$. Note that these are $\mathbb{Z}$-valued characters of $G$. 
	Under the above setup, we have following definition and results, which we use to calculate $c(G)$. 
	\begin{definition} \label{D1}
		Let $G$ be a finite group. 
		\begin{enumerate}
			\item [\rmfamily(i)] For $\psi\in \Irr(G)$, define $d(\psi)= |\Gamma(\psi)|\psi(1)$.
			\item [\rmfamily(ii)] For any complex character $\chi$ of $G$, define 
			\[  m(\chi)=
			\begin{cases}
				0 &\quad \text{ if } \chi(g) \geq 0 \text{ for all } g\in G,\\
				-\min \left\{ \sum\limits_{\sigma\in \Gamma(\chi)} \chi^{\sigma}(g): g\in G \right\} &\quad \text{ otherwise.} 
			\end{cases} \]
				\end{enumerate}
	\end{definition}

	\begin{lemma} \label{L1} Let $G$ be a non-abelian finite group. Then the following hold. 
		\begin{enumerate}
			\item [(i)] \textnormal{\cite[Lemma 2.2]{BG}} $G$ possesses a faithful quasi-permutation character of least degree, i.e. $c(G)$, which has the form $m(\xi)1_{G} + \xi$, where $\xi = \sum_{i\in I} \Psi_{i}$, $I \subseteq \{  1,2,\ldots, r \}$ and
			\begin{equation} \label{eq:lemmaL1}
				\bigcap_{i\in I} \ker(\Psi_{i}) = 1, \quad \bigcap_{i\in J} \ker(\Psi_{i}) \neq 1 \text{ for any } J \subset I.
			\end{equation}
			\item [(ii)] \textnormal{\cite[Theorem 2.3]{GA}} If $G$ is a $p$-group, then $m(\xi) = \frac{1}{p-1}\xi(1)$ and $|I|$ is the minimum number of generators of $Z(G)$.
		\end{enumerate}
	\end{lemma}
	\noindent The proof of \cite[Lemma 2.2]{BG} contains a minor error, and we correct it here.\\
	
	\noindent \emph{Proof of Lemma }\ref{L1}(i). Let $\xi$ be any faithful quasi-permutation character of $G$ of minimum degree. Then 
	\[ \xi = m'1_G + \sum_{i \in I} d_i \Psi_{i} \qquad (m', d_i \geq 1 \text{ and } I \subseteq \{ 1,\ldots, r \} ) \]
	Here, $1 = \ker(\xi) = \cap_{i\in I} \ker(\Psi_{i})$ and $-m'$ is the minimum value of $\sum_{i \in I} d_i \Psi_{i}$. Take $I' \subset I$ such that
	\begin{equation} \label{eq:BDthm2}
		\bigcap_{i\in I'} \ker(\Psi_{i}) = 1, \quad \bigcap_{i\in J} \ker(\Psi_{i}) \neq 1 \text{ for any } J \subset I'.
	\end{equation}
	(The existence of such a subset $I'$ can be proved by induction on $|I|$.) Consider the (rational valued) character $\sum_{i \in I'} \Psi_{i}$, and let $-m$ be its minimum value. Take $\eta = m1_G + \sum_{i \in I'} \Psi_{i}$. Then $\eta$ is a quasi-permutation character of $G$, and it is faithful. Hence $\eta(1) \geq \xi(1)$. We can write $\xi$ as
	\begin{equation}\label{eq:BGtheorem} 
		\xi = (m'-m)1_{G} + \sum_{i\in I', d_i >1}(d_i - 1) \Psi_{i} +   \sum_{i\in I-I'}d_{i}\Psi_{i} + \eta.
	\end{equation}
	If $m' > m$, then $\eta$ will be a sub-character of $\xi$, so $\xi(1) > \eta(1)$, a contradiction. So $m \geq m'$. Now $-m$ is the minimum value of $\sum_{i \in I'} \Psi_{i}$, say $-m = \sum_{i \in I'} \Psi_{i}(h)$ for some $h\in G$, i.e., $\eta(h) = 0$. Since $\Psi_{i}(h) \leq |\Psi_{i}(h)| \leq \Psi_{i}(1)$, we get 
	\begin{align*}
		0 \leq \xi(h) &= (m' - m) + \sum_{i\in I', d_i >1}(d_i - 1) \Psi_{i}(h) +   \sum_{i\in I-I'}d_{i}\Psi_{i}(h) \\
		&\leq (m' - m) + \sum_{i\in I', d_i >1}(d_i - 1) \Psi_{i}(1) +   \sum_{i\in I-I'}d_{i}\Psi_{i}(1)
	\end{align*}
	Therefore, from Equation \eqref{eq:BGtheorem}, we get
	\[ \xi(1) = (m'-m) + \sum_{i\in I', d_i >1}(d_i - 1) \Psi_{i}(1) +   \sum_{i\in I-I'}d_{i}\Psi_{i}(1) + \eta(1) \geq \eta(1) \geq \xi(1). \]
	So $\eta = m1_G + \sum_{i \in I'} \Psi_{i}$ is a faithful quasi-permutation character of least degree satisfying \eqref{eq:BDthm2}. \qed\\
	
	\noindent {\bf Note:} There are non-abelian groups, for example $SD_{16}$ of order 16, with a faithful quasi-permutation character of least degree, of the form $m(\xi')1_{G} + \xi'$, where $\xi' = \sum_{i\in I} d_{i}\Psi_{i}$ and $d_{i} \geq 2$ for some $i$, or $\{ \Psi_{i} \}_{i\in I}$ do not satisfy \eqref{eq:lemmaL1}.\\
		
	\noindent If $G \cong \prod_{i=1}^k C_{p^{r_i}}$ and $p\geq 3$, then $c(G)=\mu(G)=\sum_{i=1}^k p^{r_i}$ (see \textnormal{\cite[Theorem 2.11]{HB1997}} and \textnormal{\cite[Theorem 3.2]{BG}}).	
 For irreducible characters of finite $p$-groups, Ford has proved the following result in \cite{FORD}.
	
	\begin{lemma}\textnormal{\cite[Theorem 1]{FORD}} \label{thm:ford}
		Let $G$ be a $p$-group and $\chi$ an irreducible complex character of $G$. Then one of the following holds:
		\begin{enumerate}
			\item [\rmfamily(i)] There exists a linear character $\lambda$ on a subgroup $H$ of $G$ which induces $\chi$ and	generates the same field as $\chi$, i.e., $\lambda\ind_{H}^{G} = \chi$ and $\mathbb{Q}(\lambda) = \mathbb{Q}(\chi)$.
			\item [\rmfamily(ii)] $p=2$ and there exist subgroups $H < K$ in $G$ with $|K/H| = 2$ and a linear character $\lambda$ of $H$ such that with $\lambda\ind_{H}^{K} = \eta$, $[\mathbb{Q}(\lambda): \mathbb{Q}(\eta)] = 2$, $\eta\ind_{K}^{G} = \chi$, and $\mathbb{Q}(\eta) = \mathbb{Q}(\chi)$.
		\end{enumerate}
	\end{lemma}

	\section{Results} \label{sec:results}
	
	\subsection{Range of $c(G)$ for $p$-groups} \label{subsec:rangeforc(G)}
	\noindent Let us denote by $\overline{G_{p}}$, the class of $p$-groups having a minimal faithful permutation representation $\mathcal{H}=\{G_{1}, \ldots, G_{n} \}$ such that $\{ Z(G)\cap G_{1}, \ldots, Z(G)\cap G_{n} \}$ is a minimal faithful permutation representation of $Z(G)$. Johnson \textnormal{\cite[Proposition 3]{DLJ}} proved that if $G$ is a non-abelian group in $\overline{G_{p}}$, with $p\geq 3$, and $G$ is not a non-trivial direct product, then $p \mu(Z(G)) \leq \mu(G) \leq \frac{1}{p}|G:Z(G)|\mu(Z(G)).$ This result has turn out to be quite useful in the computation of $\mu(G)$ for $p$-groups.
	With the same motivation, we derive a range for $c(G)$ for any $p$-group $G$ in this subsection.\\

	\noindent	If $G$ is a non-abelian $p$-group ($p\geq 3$), then $\mu(G) = c(G) = f(p) = a_1 p + a_{2}p^2 + \cdots + a_{r}p^{r}$ with $0 \leq a_{i} < p$. In the following results, we give a sufficient condition on $G$ such that $f'(0) = 0$.

	\begin{lemma} \label{lemma:a_1=0}
		Let $G$ be a non-abelian $p$-group of order $p^{n} ~(p\geq 3)$ such that $G$ is not a direct product of an abelian and a non-abelian subgroup. Then $p^2$ divides $\mu(G)$.
	\end{lemma}
	\noindent \emph{Proof.} Let $X = \{ H_{1}, H_{2}, \ldots, H_{t} \}$ be a minimal faithful permutation representation of $G$. Then from Lemma \ref{L2}, we get $t = d(Z(G))$. Now, let $\mu(G) =\sum_{k=1}^{r} a_{k}p^{k}$, where $0 \leq a_{k} < p$ ($1\leq k \leq r$). Since $|X| = d(Z(G))$, we have $\sum_{k=1}^{r} a_{k} = d(Z(G))$.\\
	{\bf Claim: $a_{1} = 0$.} On the contrary, suppose that $a_{1} > 0$. Then for some $H_{i} \in X$, $|G/H_{i}|= p$ and $\Core_{G}(H_{i}) = H_{i}$. Without loss of generality, we can assume that $|G/H_{1}|= p$. Then, 
	\[ \bigcap_{j=1}^{t} \Core_{G}(H_{j}) = 1 \Rightarrow \Core_{G}(H_{1}) \bigcap \left( \bigcap_{j=2}^{t} \Core_{G}(H_{j})  \right) = 1 \Rightarrow H_{1} \bigcap \left( \bigcap_{j=2}^{t} \Core_{G}(H_{j})  \right) = 1.  \]
	Let $\bigcap_{j=2}^{t} \Core_{G}(H_{j}) = K_{1}$. Then $K_{1}$ is a normal subgroup of $G$ with $|K_{1}|= p$. Since $|H_{1}| = p^{n-1}$, we get $G = H_{1} \times K_{1} \cong H_{1} \times C_{p}$, a contradiction. 
	Hence, $a_{1} = 0$, which implies that $p^2$ divides $\mu(G)$ or $f'(0) = 0$. \qed 
	
	\begin{remark} \label{remark:a_1=0}
		\textnormal{ The converse of Lemma \ref{lemma:a_1=0} is not true, i.e., there exists a group $G$ such that $p^2$ divides $\mu(G)$ but $G$ is a direct product of an abelian and a non-abelian subgroup. For example, take $G =H  \times K$, where $H = \langle x, y~|~ y^{-1}x^{-1}yx = x^{p}, x^{p^2} = y^p = 1 \rangle$ and $K \cong C_{p^2}$. Here, $H$ is a non-abelian $p$-group of order $p^3$. Then from \cite[Corollary 2.2]{DW}, $\mu(G) = \mu(H) + \mu(K) = 2p^2$, and hence $p^2$ divides $\mu(G)$.
		}
	\end{remark}

		\noindent {\bf Proof of Theorem \ref{lemmalinearcharacter}.} Suppose $X_G = \{ \chi_{i} \}_{i=1}^{m}$, where $\chi_{i} \in \nl(G)$, for each $1\leq i \leq m$. Set for each $1\leq i \leq m$, $\chi_{i}(1) = p^{a_{i}}$ and $|\Gamma(\chi_{i})| = \phi(p^{b_i})$. Then
		\[ \xi_{X_G}(1) = \sum_{i=1}^{m} \sum_{\sigma \in \Gamma(\chi_{i})} \chi_{i}^{\sigma} (1) = \sum_{i = 1}^{m} \chi_{i}(1)|\Gamma(\chi_{i})| = \sum_{i = 1}^{m} p^{a_{i}} \phi(p^{b_{i}}) = \sum_{i = 1}^{m} p^{a_{i}+b_{i}} - p^{a_{i}+b_{i} - 1} . \]
		From Lemma 12, 
		\[ m(\xi_{X_G}) = \frac{1}{(p-1)} \xi_{X_G}(1) = \sum_{i = 1}^{m} p^{a_{i}+b_{i} - 1}. \]
		Thus, we get $c(G) = \xi_{X_G}(1) + m(\xi_{X_G}) = \sum_{i = 1}^{m} p^{a_{i}+b_{i}}.$
		Now, suppose $p^s = \min\{ \chi(1) ~|~ \chi \in \nl(G)) \}$. Since $a_{i} \geq s$ and $b_{i} \geq 1$, for each $1\leq i \leq m$, we get $p^{s+1}$ divides $c(G)$.\\
		Now, we prove that if $d(Z(G)  \cap G') = d(Z(G))$, then $X_{G} \cap \lin(G) = \emptyset$.
		When $d(Z(G)) = 1$, i.e., $Z(G)$ is cyclic, then the result is obvious. Now, suppose $d(Z(G)) = m > 1$ with $d(Z(G)  \cap G') = d(Z(G))$. Let $X_G = \{ \psi_{i} \}_{i=1}^{m}$ be a minimal faithful quasi-permutation representation of $G$ satisfying \eqref{eq:X_G}. To show that each $\psi_{i}$ is nonlinear, suppose $\psi_{r}$ is linear, for some $1\leq r \leq m$. Since $\cap_{i=1, i\neq r}^{m} \ker(\psi_{i}) \neq 1$, take $z \in \left(\cap_{i=1, i\neq r}^{m} \ker(\psi_{i})\right) \cap Z(G)$ of order $p$. Since $d(Z(G) \cap G') = d(Z(G))$, so all the central elements of order $p$ in $Z(G)$ are also in $G'$. Hence, $z\in G' \subseteq \ker(\psi_{r})$, and so $z \in \cap_{i=1}^{m} \ker(\psi_{i})$, which is a contradiction. \qed

		\begin{remark} \label{remark:NotInX_G}
			\textnormal{ In Theorem \ref{lemmalinearcharacter}, if $d(Z(G) \cap G') \neq d(Z(G))$, then $\psi$ may belong in $X_{G}$, for some $\psi\in \lin(G)$. For example, consider
				\begin{align*}
					G = \langle \alpha_{1}, \alpha_{2}, \alpha_{3}, \alpha_{4}, \alpha_{5}, \beta_{1}, \beta_{2} ~|~ & [\alpha_{4}, \alpha_{5}] = \alpha_{3}, [\alpha_{3}, \alpha_{5}] = \alpha_{2}, [\alpha_{2}, \alpha_{5}] = \alpha_{1} = \beta_{1}, \alpha_{4}^{p}= \beta_{2},\\
					& \alpha_{5}^{p} = \beta_{1}, \alpha_{2}^{p} = \alpha_{3}^{p} =   \beta_{1}^{p} = \beta_{2}^{p} =  1 \rangle
				\end{align*}
				which is a $p$-group of order $p^6$ ($p\geq 5$) belonging to the isoclinic family $\Phi_{9}$ (see \cite{Easterfield}). In the above group presentation, all relations of the form $[\alpha, \beta] = 1$ (with $\alpha$, $\beta$ generators) have been omitted.
				Now, here $Z(G) = \langle \alpha_{5}^{p}, \alpha_{4}^{p} \rangle \cong C_{p} \times C_{p}$, $G{}' = \langle \alpha_{5}^{p}, \alpha_{3}, \alpha_{2} \rangle \cong C_{p}^{3}$ and $G/G{}' = \langle \alpha_{5} G{}', \alpha_{4} G{}' \rangle \cong C_{p} \times C_{p^2}$. Hence $d(Z(G) \cap G') = 1 \neq d(Z(G))$. Now, take $A = \langle \alpha_{5}^{p}, \alpha_{4}, \alpha_{3}, \alpha_{2} \rangle \cong C_{p} \times C_{p^2} \times C_{p} \times C_{p}$. Then $A$ is an abelian normal subgroup of index $p$ in $G$. Hence, from \cite[Theorem 12.11]{I}, $\cd(G) = \{ 1, p \}$. From Lemma \ref{L1}, $|X_{G}| = 2$.  \\
				Let $\psi\in\lin(G/G{}')$ given by $\psi =1_{\langle \alpha_{5} G{}' \rangle} \cdot \psi_{\langle \alpha_{4}G{}' \rangle}$, where $\psi_{\langle \alpha_{4} G{}' \rangle}$ is a faithful linear character of $\langle \alpha_{4} G{}' \rangle$. On the other hand, let $\lambda \in \lin(A)$ given by $\lambda = \lambda_{\langle \alpha_{5}^{p} \rangle} \cdot 1_{\langle \alpha_{4} \rangle} \cdot 1_{\langle \alpha_{3} \rangle} \cdot 1_{\langle \alpha_{2} \rangle}$, where $\lambda_{\langle \alpha_{5}^{p} \rangle}$ is a faithful linear character of $\langle \alpha_{5}^{p} \rangle$. Since the inertia group of $\lambda$ in $G$, $\inertiagroup_{G}(\lambda) = A$, we get $\lambda\ind_{A}^{G} \in \nl(G)$. Here, $\ker(\psi) \cap \ker(\lambda\ind_{A}^{G}) = \ker(\psi) \cap \Core_{G}(\ker(\lambda)) = 1$. Now, $d(\psi) = \phi(p^2) = p(p-1)$, and $d(\lambda\ind_{A}^{G}) = \lambda\ind_{A}^{G}(1) |\Gamma(\lambda\ind_{A}^{G})| \leq p|\Gamma(\lambda)| = p\phi(p) = p(p-1)$. Since $|\Gamma(\lambda\ind_{A}^{G})| \geq \phi(p)$, we get $d(\lambda\ind_{A}^{G}) = p(p-1)$. Suppose $ \xi =  \left[ \sum_{\sigma \in \Gamma(\psi)} \psi^{\sigma}  \right] + \left[ \sum_{\sigma \in \Gamma(\lambda\ind_{A}^{G})} \left( \lambda\ind_{A}^{G} \right)^{\sigma}  \right]$. Then $\xi(1) = 2p(p-1)$. From \cite[Lemma 4.5]{HB}, we get $m(\xi) = 2p$. Hence, $\xi(1) + m(\xi) = 2p^2$, which implies that $c(G) \leq 2p^2$. 
				From Lemma \ref{lemma:a_1=0} and \cite[Theorem 3.2]{BG}, it is easy to see that $c(G) = \mu(G) \geq 2p^2$. Therefore, we get $c(G) = 2p^2$ and $X_{G} = \{ \psi, \lambda\ind_{A}^{G} \}$, where $\psi \in \lin(G)$. }
		\end{remark}
		
		\begin{remark} \label{remark:X_Gcontainslinear}
			\textnormal{There are non-abelian $p$-groups in which there is no minimal faithful quasi-permutation representation $X_G$ such that $X_G \cap \lin(G) = \emptyset$. For example, consider $G = H \times K$, where $H$ is an extraspecial $p$-group of order $p^3$ ($p\geq 3$) with $\exp(H) = p$, and $K\cong C_p$. From \cite[Theorem 3.2]{BG}, $c(G) = p^2+p$. If possible, suppose $X_G = \{ \chi_{1}, \chi_{2} \} \subset \nl(G)$ be a minimal faithful quasi-permutation representation of $G$. Then \[ \xi(1) = \sum_{i=1}^{2} \sum_{\sigma \in \Gamma(\chi_i)}{\chi_{i}}^{\sigma}(1) = 2p(p-1) \text{ and } m(\xi) = 2p \text{ (from Lemma \ref{L1})}, \]
				which implies that $2p^2 = \xi(1) + m(\xi) > c(G)$, a contradiction. }
		\end{remark}
		
		\begin{corollary} \label{cor:NotInX_G2}
			Let $G$ be a non-abelian $p$-group with $\exp(G) = p$, $\cd(G) = \{ 1, p^s \} ~ (s>1)$ and $d(Z(G) \cap G') = d(Z(G))$. Then $c(G) = d(Z(G)) p^{s+1}$.
		\end{corollary}
		
		\noindent \emph{Proof.} From Lemma \ref{L1}, we get $|X_{G}|=d(Z(G)) = t$ (say). Let $X_{G} = \{ \chi_{i} \}_{i=1}^{t}$. From Theorem \ref{lemmalinearcharacter}, we get $X_{G} \cap \lin(G) = \emptyset$. Hence, we get $\chi_{i}(1) = p^s$ for $1\leq i \leq t$. Since $G$ is an $M$-group, for $1\leq i \leq t$, there exists $H_{i} \leq G$ with $|G/H_{i}| = p^s$ such that $\chi_{i} = \lambda_{i}\ind_{H_{i}}^{G}$, for some $\lambda_{i} \in \lin(H_{i})$. Since $\exp(G) = p$ and $\lambda_{i}$ is not a trivial character of $H_i$, we get $\mathbb{Q}(\chi_{i}) = \mathbb{Q}(\lambda_{i}) = \mathbb{Q}(\omega_{p})$, for $1\leq i \leq t$. This implies that $d(\chi_{i}) = \chi_{i}(1) |\Gamma(\chi_{i})| = p^s |\Gamma(\lambda_{i})| = p^s \phi(p) = p^{s+1} - p^{s}$. From Lemma \ref{L1}, we get $c(G) = t p^{s+1} = d(Z(G)) p^{s+1}$. \qed 
		
		\begin{proposition}	\label{cor:NotInX_G}
			Let $G$ be a non-abelian $p$-group such that $G$ has an elementary abelian subgroup of index $p$ and $X_{G} \cap \lin(G) = \emptyset$. Then $c(G) = d(Z(G))p^2$.
		\end{proposition}
		\noindent \emph{Proof.} From \cite[Theorem 12.11]{I}, we get $\cd(G) = \{ 1, p \}$. By Lemma \ref{L1}, we get $|X_{G}|=d(Z(G))$. Since $X_{G} \cap \lin(G) = \emptyset$, we get $\chi(1) = p$ for all $\chi \in X_{G}$. Let $A$ be an elementary abelian normal subgroup of index $p$ in $G$. Take any $\chi \in X_G$. If $\lambda \in \Irr(A)$ appears in $\chi\restr_{A}$, then $\chi$ appears in $\lambda\ind_{A}^{G}$. Since $\lambda\ind_{A}^{G}(1) = p$, we get that $\lambda\ind_{A}^{G} = \chi$. Hence $\mathbb{Q}(\chi) = \mathbb{Q}(\lambda) = \mathbb{Q}(\omega_{p})$ and $|\Gamma(\chi)| = p-1$. Therefore $d(\chi) = p(p-1)$ for all $\chi \in X_G$.
		From Lemma \ref{L1}, it is easy to see that $c(G) = d(Z(G)) p^2$. \qed \\
			
		\noindent Now we prove Theorem \ref{P2}.\\
		
		\noindent {\bf Proof of Theorem \ref{P2}.} 
		By Lemma \ref{L1}, we get $|X_{G}|= d(Z(G)) = m$ (say). Now, let $X_{G} = \{ \chi_{1}, \chi_{2},\ldots,\chi_{m} \}$ and suppose $\xi =\sum_{i=1}^m (\sum_{\sigma\in \Gamma(\chi_{i})}\chi_{i}^{\sigma})$. Since $\exp(G) = p^b$, so $\chi_{i}(g) \in \mathbb{Q}(\omega_{p^b})$ for all $1\leq i \leq m$ and all $g\in G$. Hence $|\Gamma(\chi_{i})| \leq \phi(p^b) = p^{b} - p^{b-1}$.
		Further, $\chi_{i}(1) \leq p^e$ for all $i$. Hence
		\[ c(G) = \xi(1) + m(\xi) =\left( \sum_{i=1}^{m} \sum_{\sigma\in \Gamma(\chi_{i})}\chi_{i}^{\sigma} (1) \right) + \frac{1}{p-1} \left(\sum_{i=1}^{m} \sum_{\sigma\in \Gamma(\chi_{i})}\chi_{i}^{\sigma} (1)\right) \leq \frac{p}{p-1} \sum_{i=1}^{m} p^e(p^b - p^{b-1}) = mp^{b+e}.   \]
		On the other hand, if $g\in G$ is of order $p^b$, then $p^b = c(\langle g \rangle ) \leq c(G)$. \\
		Now, let $d(\chi_{i}) = |\Gamma(\chi_{i})|\chi_{i}(1) = \phi(p^{t_{i}})$, for $1\leq i \leq m$ and $1\leq t_{i} \leq b+e$. Then it is easy to see that $\xi(1)= \sum_{i=1}^{m} p^{t_{i}} - p^{t_{i}-1}$. From Lemma \ref{L1}, we get $m(\xi) = \sum_{i=1}^{m} p^{t_{i}-1}$. Hence, we get $c(G) = \sum_{i=1}^{m} p^{t_{i}}$, for $1\leq t_{i} \leq b+e$. This implies that $c(G) = f(p) = a_1 p + \cdots + a_{r}p^{r}$, where $0 \leq a_{i} < p$ for each $i$, and $a_r \neq 0$ with $\sum_{l = 1}^{r} a_{l} = m$.
		Then
		\begin{enumerate}
			\item [(i)] $p^b \leq c(G) =  a_1 p + a_2 p^2 + \cdots + a_{r}p^{r}$ implies $b \leq r$.
			\item [(ii)] $c(G) \leq mp^{b+e} $ (i.e. $a_1 p + a_2 p^2 + \cdots + a_{r}p^{r} \leq mp^{b+e})$ implies if $m < p$ then $r \leq b+e$ (for, if $b+e < r$ then $p^{b+e+1} \leq p^r \leq a_1 p + a_2 p^2 + \cdots + a_{r}p^{r} < p^{b+e+1}$, a contradiction).
		\end{enumerate}

		Now to prove that the range $b \leq r \leq b+e$ is best possible, we consider the following groups for $p\geq 5$:
		\begin{align*}
			G_{1} = \langle \alpha_{1}, \ldots, \alpha_{5},  \beta_{1} ~|~ [\alpha_{4}, \alpha_{5}] =& \alpha_{3}, [\alpha_{3}, \alpha_{5}] = \alpha_{2}, [\alpha_{2}, \alpha_{5}] = [\alpha_{3}, \alpha_{4}] = \alpha_{1} = \beta_{1}^{p},\\
			& \alpha_{2}^{p} = \alpha_{3}^{p} =  \alpha_{4}^{p} = \alpha_{5}^{p}= \beta_{1}^{p^2} =  1 \rangle \\
			G_{2} =	\langle \alpha_{1}, \ldots, \alpha_{5},  \beta_{1}, \beta_{2} ~|~ [\alpha_{4}, \alpha_{5}] =& \alpha_{3}, [\alpha_{3}, \alpha_{5}] = \alpha_{2}, [\alpha_{2}, \alpha_{5}] = [\alpha_{3}, \alpha_{4}] = \alpha_{1} = \beta_{1}, \alpha_{5}^{p}= \beta_{2},\\
			& \alpha_{2}^{p} = \alpha_{3}^{p} =  \alpha_{4}^{p} = \beta_{1}^{p} = \beta_{2}^{p} =  1 \rangle \text{ and }\\
			G_{3} =	\langle \alpha_{1}, \cdots, \alpha_{6} ~|~ [\alpha_{3}, \alpha_{4}] =& \alpha_{1}, [\alpha_{5}, \alpha_{6}] = \alpha_{2}, \alpha_{3}^{p} = \alpha_{1}, \alpha_{1}^{p} = \alpha_{2}^{p} =  \alpha_{4}^{p} = \alpha_{5}^{p}= \alpha_{6}^{p} =  1 \rangle.
		\end{align*}
		The groups $G_{1}, G_{2}$ and $G_{3}$ mentioned above are of order $p^6$. Here, $G_{1}$ and $G_{2}$ belong to the isoclinic family $\Phi_{10}$, whereas $G_{3}$ belongs to the isoclinic family $\Phi_{12}$ (see \cite{Easterfield}). In all the group presentations, all relations of the form $[\alpha, \beta] = 1$ (with $\alpha$, $\beta$ generators) have been omitted. Here $\exp(G_{i}) = p^2$ and $\cd(G_{i}) = \{ 1,p,p^2 \}$, for $1\leq i \leq 3$. Then $b = 2$ and $e = 2$ for all the three groups. We have the following cases.
		\begin{enumerate}
			\item $\mathbf{r = b+e:}$ 
			We have $Z(G_{1}) = \langle \beta_{1} \rangle \cong C_{p^{2}}$. By proceeding along the lines of Remark \ref{remark:NotInX_G}, we get $c(G_{1}) = p^{4}$. Then, $r = 4$, $b=2$, and  $e=2$. Hence, we get $r= 4 = 2+2 = b+e$. 
			\item $\mathbf{b < r < b+e:}$ 
			We have $Z(G_{2}) = \langle \beta_{1}, \beta_{2} \rangle \cong C_{p} \times C_{p}$. By proceeding along the lines of Remark \ref{remark:NotInX_G}, we get $c(G_{2}) = p^{3}+p^{2}$. Then, $r =3$, $b=2$, and  $e=2$. Hence, we get $b< r <  b+e$.
			\item $\mathbf{b= r:}$
			We have $Z(G_{3}) = \langle \alpha_{1}, \alpha_{2} \rangle \cong C_{p} \times C_{p}$. By proceeding along the lines of Remark \ref{remark:NotInX_G}, we get $c(G_{3}) = 2p^{2}$. Then, $r = b =2$.
		\end{enumerate}
		This concludes the proof.  \qed 
		
		\begin{remark} \label{remark:rangec(G)}
			\textnormal{Theorem \ref{P2} can also be expressed in the following form: \\
				Let $G$ be a non-abelian $p$-group such that $d(Z(G)) < p$, $\exp(G) = p^{b}$, and $\max\cd(G) = p^{e}$. Then there exist $f(p) = \sum_{k=1}^{b}a_{k}p^{k}$ and $g(p) = \sum_{l=1}^{b+e}c_{l}p^{l}$ with $0 \leq a_{k}, c_{l} < p$,  $a_{b}\neq 0, c_{b+e}\neq 0$ and $ \sum_{k=1}^{b}a_{k} = \sum_{l=1}^{b+e}c_{l} = d(Z(G))$ such that $f(p) \leq c(G) \leq g(p)$. } 
				\end{remark}

		\noindent {\bf Proof of Corollary \ref{cor:newrange2}.} From Theorem \ref{P2} and Lemma \ref{lemma:a_1=0}, the result follows. \qed

		\subsection{Finite non-abelian $p$-groups with cyclic center} \label{subsec:cycliccenter}
		In Theorem \ref{thm:cyclicexistence}, we prove the existence of a non-abelian permutation group of order $p^n$ acting transitively on $p^{i}$ points, for each $i$ ($2 \leq i \leq n-1$).\\
		
		\noindent	{\bf Proof of Theorem \ref{thm:cyclicexistence}.}
		Let $G$ be a non-abelian $p$-group of order $p^{n}$ ($p\geq 3$) with cyclic center. Since $d(Z(G)) = 1$, from Lemma \ref{L2}, all the minimal faithful permutation representations of $G$ are transitive. Hence, we get $p^2 \leq \mu(G) \leq p^{n-1}$. Therefore, $\mu(G) = c(G) \in \{ p^2, p^3, \ldots, p^{n-1} \}$.\\ Now, suppose $p\geq n$. We show the existence of a group $G$ of order $p^n$ such that $c(G) = \mu(G) = p^i$, for each $2 \leq i \leq n-1$. 
		In the following group presentations, all relations of the form $[\alpha, \beta] = 1$ (with $\alpha$, $\beta$ generators) have been omitted.  The symbol $\alpha_{i+1}^{(p)}$ denotes $\alpha_{i+1}^{p}\alpha_{i+2}^{\binom{p}{2}}\cdots\alpha_{i+k}^{\binom{p}{k}} \cdots \alpha_{i+p}$, where $i$ is a positive integer and $\alpha_{i+2}, \ldots,\alpha_{i+p}$ are suitably defined. Consider
		\[ G_{2} = \langle \alpha, \alpha_{1}, \ldots, \alpha_{n-1} ~|~ [\alpha_{j}, \alpha]= \alpha_{j+1}, \alpha^{p} = \alpha_{j}^{(p)} = \alpha_{n-1}^{(p)} = 1~  (j=1, 2,\ldots, n-2) \rangle. \]
		Since $p\geq n$, it is easy to check that $\alpha_{j}^{(p)} = \alpha_{j}^{p} = 1$ ($1\leq j \leq n-1$). Here, $|G_{2}| = p^{n}$, $\cd(G_{2}) = \{ 1, p \}$, $Z(G_{2}) = \langle \alpha_{n-1} \rangle \cong C_{p}$ and the nilpotency class of $G_{2}$ is $n-1$.
		Now, take $H = \langle \alpha_{1},\alpha_{2},\ldots,\alpha_{n-2} \rangle \cong C_{p}^{n-2}$. Since $H \cap Z(G_{2}) = 1$, we get $\Core_{G_{2}}(H) = 1$. Thus, $\mu(G_{2}) \leq p^2$. Since $G_{2}$ is a non-abelian group, we get $\mu(G_{2}) \geq p^2$. Therefore, $\mu(G_{2}) = p^2$.\\
		Now, for $3\leq i \leq n-1$, consider
		\begin{align*}
			G_{i} = \langle \alpha, \alpha_{1}, \alpha_{2}, \ldots, \alpha_{n-i+1}, \alpha_{n-i+2} ~|~ [\alpha_{1}, \alpha]=& \alpha_{2}, [\alpha_{j}, \alpha]= \alpha_{j+1}, \alpha_{1}^{p^{i-2}} = \alpha_{n-i+2},\\
			& \alpha^{p} = \alpha_{j}^{(p)} = \alpha_{n-i+2}^{(p)} = 1, ~ ( j= 2, 3, \ldots, n-i+1) \rangle.
		\end{align*}
		It is easy to check that for $p\geq n-1$, we have $\alpha_{j}^{(p)} = \alpha_{j}^{p} = 1$ ($2 \leq j \leq n-i+2$). Here, for $3\leq i \leq n-1$, $|G_{i}| = p^n$, $\cd(G_{i}) = \{ 1, p \}$, $Z(G_{i}) = \langle \alpha_{1}^{p} \rangle \cong C_{p^{i-2}}$ and the nilpotency class of $G_{i}$ is $n+2 - i$.
		Now, take $H = \langle \alpha_{2},\alpha_{3},\ldots,\alpha_{n- i + 1} \rangle \cong C_{p}^{n-i}$.
		Since $H \cap Z(G) = 1$, we get $\Core_{G_{i}}(H) = 1$. Thus, $\mu(G_{i}) \leq p^{n-(n-i)} = p^{i}$, for each $i$ ($3\leq i \leq n-1$). On the other hand, $\mu( \langle \alpha_{1}, \alpha_{2} \rangle ) = p^{i-1}+p $, so $\mu(G_{i}) \geq p^{i}$, for $3\leq i \leq n-1$. Therefore, we get $\mu(G_{i}) = p^{i}$, for each $i$ ($3\leq i \leq n-1$). \qed

		\begin{lemma} \label{lemma:cdG3}
			Let $G$ be a non-abelian $p$-group such that  $\cd(G) = \{ 1, p, p^a \}$, for some integer $a>1$. If $\chi \in \Irr(G)$ is faithful, then $\chi(1) = p^a$.
		\end{lemma}
		\noindent \emph{Proof.} Let $\chi$ be a faithful irreducible character of $G$, and suppose $\chi(1) = p$. Now $\chi = \lambda\ind_{H}^{G}$ for some subgroup $H$ of index $p$ (since $G$ is an $M$-group) and for some $\lambda \in \lin(H)$. Then $H'$ is normal in $G$ and is contained in $\ker(\lambda)$. So $H' \subset \ker(\chi) = 1$, i.e., $H$ is abelian of index $p$ in $G$. Then $\cd(G) \subset \{ 1, p \}$, which is a contradiction. \qed \\
		
		\noindent Note that Isaacs and Moret\'{o} \cite[Theorem 5.2]{IM} proved the existence of a $p$-group with character degree set $\{ 1, p, p^a \}$, for any integer $a>1$ and a prime $p$.\\

		\noindent Now, with the help of Theorem \ref{P2}, we prove Corollary \ref{cor:newrangecyclic}, in which we compute the minimal faithful quasi-permutation representation degree for a general $p$-group with cyclic center.\\
		
		\noindent {\bf Proof of Corollary \ref{cor:newrangecyclic}.} From Lemma \ref{L1}, we get $|X_{G}| = 1$. Suppose $Z(G) \cong C_{p^m}$, $p^{e} = \max \cd(G)$ and $p^{\alpha} = \min\{ \chi(1)~|~ \ker(\chi) = 1, \chi \in \nl(G) \}$. Let $X_G = \{ \chi \}$, where $\chi$ is a faithful irreducible character of $G$. Since $Z(G)$ is cyclic, so $d(\chi) = \chi(1) |\Gamma(\chi)| \geq p^{\alpha} \phi(|Z(G)|) = p^{\alpha}\phi(p^m) = p^{\alpha +m-1}(p-1)$. From Lemma \ref{L1}, we get $p^{\alpha+m} = p^{\alpha} |Z(G)| \leq c(G)$. From Theorem \ref{P2}, we get $c(G) \leq p^{e} \cdot \exp(G)$. Therefore, $p^{\alpha} |Z(G)| \text{ divides } c(G)$ and $c(G) \text{ divides } p^{e} \cdot \exp(G)$. 
		Now, suppose $\cd(G) = \{ 1, p, p^{e} \}$ ($e>1$) and $\exp(G) = p$. 
		Then from Lemma \ref{lemma:cdG3}, $\chi(1) = p^e$. By the above discussion, we get $p^{e} |Z(G)| \leq c(G) \leq p^{e} \cdot p$. This implies that $c(G) = p^{e+1}$. \qed \\
		
		\noindent We now compute $c(G)$ for a normally monomial $p$-group $G$ with cyclic center. Suppose $\chi$ is a faithful irreducible character of $G$. Then all the normal subgroups which linearly induce $\chi$ are abelian and have maximal order among all abelian subgroups of $G$ (see \cite[Proposition 3]{AM}). Moreover, all faithful irreducible characters of $G$ have the same degree, which is the maximal degree of all irreducible characters of $G$. Note that metabelian groups are normally monomial groups (see \cite{BGB}). Now we prove Theorem \ref{lemma:normallymonomial}.\\
		
		\noindent {\bf Proof of Theorem \ref{lemma:normallymonomial}:}
		Let $X_G = \{ \chi\}$, where $\chi$ is a faithful irreducible character of $G$. Suppose $A$ is an abelian normal subgroup of maximum order in $G$. Then $\chi = \lambda\ind_{A}^{G}$, for some $\lambda \in \lin(A)$, and $\chi(1) = \max\cd(G)$ \cite[Proposition 3]{AM}. Let $|\Gamma(\chi)| = \phi(p^b)$, for some $b\geq 1$. Then $d(\chi) = \chi(1) |\Gamma(\chi)| = \chi(1)\phi(p^b)$. From Lemma \ref{L1}, we get $c(G) = \chi(1)p^b$. Since $\mathbb{Q}(\chi) \subseteq \mathbb{Q}(\lambda) \subseteq \mathbb{Q}(\omega_{\exp(A)})$, we get $p^b \leq \exp(A)$. Further, since $Z(G)$ is cyclic and $\chi$ is faithful, $\phi(|Z(G)|) \leq |\Gamma(\chi)| = \phi(p^b)$. Thus, $|Z(G)| \leq p^b$. Therefore, we get $\chi(1)|Z(G)| \leq c(G) = \chi(1)p^b \leq \chi(1)\exp(A)$, and hence $(\max\cd(G))|Z(G)| \leq c(G) \leq (\max\cd(G))\exp(A)$. This completes the proof. \qed

		\begin{remark} \label{remark:rewritenormMon}
			\textnormal{ Theorem \ref{lemma:normallymonomial} can also be expressed as follows:\\
				Let $G$ be a normally monomial $p$-group of order $p^n$ with $\max \cd(G) = p^{a}$, and $Z(G) \cong C_{p^{m}}$. Suppose $A$ is an abelian normal subgroup of maximum order in $G$. If $\exp(A) = b$, then $c(G) \in \{ p^{a+m}, p^{a+m+1}, \ldots, p^{a+b} \}$. }
		\end{remark}

		\begin{corollary} \label{cor:normonomial2} \label{cor:normallymonomial}
			Let $G$ be a normally monomial $p$-group with cyclic center.
			\begin{enumerate}
				\item [(i)] If $\exp(G) = p$ and $\max \cd(G) = p^{a}$, then $c(G) = p^{a+1}$.
				\item [(ii)] If $A$ is an abelian normal subgroup of maximum order in $G$ such that $\exp(A) = |Z(G)|$, then
				$c(G) = (\max \cd(G))|Z(G)|$.
			\end{enumerate}
			
		\end{corollary}
		\noindent \emph{Proof.} From Remark \ref{remark:rewritenormMon}, (i) is immediate.
		From Theorem \ref{lemma:normallymonomial}, (ii) follows. \qed\\
				
		 If $G$ is a finite $p$-group ($p\geq 3$) of nilpotency class 2, then $G$ is a normally monomial $p$-group. In this case, we derive the following result.
				
		\begin{proposition}\label{lemma:class2}
			Let $G$ be a finite $p$-group $(p\geq 3)$ of nilpotency class 2 with cyclic center. Then there exists an abelian normal subgroup, say $H$, of maximum order in $G$ such that $\exp(H) = |Z(G)|$.
		\end{proposition}
		\noindent \emph{Proof.} Let $G$ be a finite $p$-group ($p\geq 3$) of nilpotency class 2 with cyclic center. Let $\chi$ be a faithful irreducible character of $G$.
		From Lemma \ref{thm:ford}, there exists a linear character $\lambda$ on a subgroup $H$ of $G$ such that $\chi = \lambda\ind_{H}^{G}$ and $\mathbb{Q}(\chi) = \mathbb{Q}(\lambda)$. Here, $H{}' \subseteq \ker(\lambda)$. Since $G$ is a group of nilpotency class 2, we get $H' \subseteq G' \subseteq Z(G)$. Since $\ker(\chi) = \Core_{G}(\ker(\lambda)) = 1$, $H{}' = 1$.	Since $G$ is a metabelian $p$-group, $G$ is a normally monomial $p$-group. Let $A$ be a normal subgroup which linearly induces $\chi$. Then from \cite[Proposition 3]{AM}, $A$ is an abelian normal subgroup of  maximum order in $G$. Moreover, $\chi(1) = \max \cd(G)$. Since $\chi(1) = |G/A| = |G/H|$, $H$ also has maximum
		order among all abelian subgroups of G. Hence, $G{}' \subseteq Z(G) \subseteq H$, which implies that $H$ is an abelian normal subgroup of $G$.
		
		Suppose $X_G = \{ \chi \}$, where $\chi$ is a faithful irreducible character of $G$. Then from the above discussions, there exists an abelian normal subgroup $H$ of maximum order in $G$ and $\lambda \in \lin(H)$ such that 
		$\chi = \lambda\ind_{H}^{G} \text{ and } \mathbb{Q}(\chi) = \mathbb{Q}(\lambda)$.
			Since $G$ is a metabelian group, we get $\left| \frac{H}{\ker(\lambda)} \right| = \exp(H)$. This implies that $|\Gamma(\chi)| = |\Gamma(\lambda)| = \phi(\exp(H))$.
%
	 On the other hand, 		
		since $G$ is of nilpotency class 2, $\chi(1) = |G/Z(G)|^{1/2}$ (see \cite[Theorem 2.31]{I}), and $\chi(g) = 0$ for all $g\in G\setminus Z(G)$ (see \cite[Corollary 2.30]{I}). Since $Z(G)$ is cyclic, $\chi\restr_{Z(G)} = \chi(1)\mu$ for some faithful linear character $\mu$ of $Z(G)$. Hence, $|\Gamma(\chi)| = |\Gamma(\mu)| = \phi(|Z(G)|)$. Therefore, $\phi(|Z(G)|) = \phi(\exp(H))$. This completes the proof.	  \qed

		\begin{consequence} \label{remark:behraveshresult}
			\begin{enumerate}
				\item \textnormal{If $G$ is a finite $p$-group ($p\geq 3$) of nilpotency class $2$ with cyclic center, then by Theorem \ref{lemma:normallymonomial} and Proposition \ref{lemma:class2}, it follows that $c(G)=(\max\cd(G))|Z(G)|=|G/Z(G)|^{1/2}|Z(G)|$. This shows that Theorem \ref{thm_HB_class2} can be deduced from Theorem \ref{lemma:normallymonomial} for $p\geq 3$. } 
				\item \textnormal{In the special case, when $G$ is a normally monomial $p$-group with cyclic center and $G'$ is an abelian normal subgroup of maximum order in $G$, then $G'$ is the unique normal abelian normal subgroup of maximum order in $G$. Thus, in this case, we get} 
				\[ (\max \cd(G))|Z(G)| \text{ divides } c(G) \text{ and } c(G) \text{ divides } (\max \cd(G))\exp(G').\]
			\end{enumerate}
		\end{consequence}

		\begin{corollary} \label{cor:exp(G')}
			Let $G$ be a non-abelian $p$-group of order $p^n$ $(p\geq 5 \text{ and } n \leq 6)$ with cyclic center, and suppose $G'$ is abelian with maximal order among all abelian subgroups of $G$. Then the following hold.
			\begin{enumerate}
				\item $n>4$.
				\item When $n=5$, $c(G) = p^3$ and when $n = 6$, $c(G) \in \{ p^3, p^4 \}$.
			\end{enumerate}
		\end{corollary}
		\noindent \emph{Proof.} It is well known that $p$-groups of order $\leq p^4$ contain an abelian subgroup of index $p$ (see \cite[Exercise 8, p. 28]{YBbook}). Also for a non-abelian $p$-group $G$, we have $|G:G'| \geq p^2$. Hence $n> 4$.\\
		Now, let $|G| = p^5$ ($p\geq 5$). Then from \cite[Section 4.5]{RJ}, it is easy to see that if $G\in \Phi_{i}$, where $i\in \{ 2,3,\ldots,9 \}$, then $G'$ does not have maximal order among all abelian subgroups of $G$. Now, let $G\in \Phi_{10}$. Then $G'\cong C_{p}^{3}$, $\cd(G) = \{ 1, p, p^2 \}$ and $Z(G) \cong C_{p}$ (see \cite[Sections 4.1 and 4.5]{RJ}). Thus, $\exp(G') = p =|Z(G)|$, and hence from Theorem \ref{lemma:normallymonomial}, $c(G) = \max\cd(G) \cdot \exp(A) = p^2\cdot p = p^3$.\\
		Now, let $|G| = p^6$ ($p\geq 5$). Then from \cite{Easterfield}, it is easy to see that if $G\in \Phi_{i}$, where $i\in \{ 2,3,\ldots,43 \} \setminus \{ 36,38,40,41,42,43 \}$, then $G'$ does not have maximal order among all abelian subgroups of $G$. Now, let $G\in \Phi_{i}$, where $i\in \{ 36,38,40,41 \}$. Then $G'\cong C_{p}^{4}$, $\cd(G) = \{ 1, p, p^2 \}$ and $Z(G) \cong C_{p}$ (see \cite{Easterfield}). Thus, $\exp(G') = p =|Z(G)|$, and hence from Theorem \ref{lemma:normallymonomial}, we get
		\[ c(G) = \max\cd(G) \cdot \exp(A) = p^2\cdot p = p^3, \text{ for all } G\in \Phi_{i}, \text{ where } i\in \{ 36,38,40,41 \} \text{ and } |G| = p^6 (p\geq 5). \]
		Next, suppose $G \in \Phi_{i}$, where $i\in \{42, 43 \}$. Then $G{}' \cong C_{p^2} \times C_{p} \times C_{p}$, $\cd(G) = \{ 1, p, p^2 \}$ and $Z(G) \cong C_{p}$ (see \cite{Easterfield}). Here $\exp(G') = p^2$. Thus, from Theorem \ref{lemma:normallymonomial}, we get 
		\[ (\max \cd(G))|Z(G)| \leq c(G) \leq (\max \cd(G))\exp(A) \Rightarrow p^3 \leq c(G) \leq p^4. \] 
		To get exact value of $c(G)$, we need to study each group from $\Phi_{42}$ and $\Phi_{43}$ separately. We use group presentations computed by Easterfield in \cite{Easterfield}.  The notation $G_{(i,j)}$ denotes that the group is the $j^{th}$ group in the isoclinic family $\Phi_{i}$ in the list of groups of order $p^6$ ($p\geq 5$). All relations of the form $[\alpha, \beta] = 1$ (with $\alpha$, $\beta$ generators) have been omitted in the following group presentations. We have the following groups.
		\begin{align*}
			G_{(42,1)} = &  \langle \alpha_{1}, \alpha_{2}, \alpha_{3}, \alpha_{4}, \alpha_{5}, \alpha_{6} ~|~ [\alpha_{5}, \alpha_{6}] = \alpha_{4}, [\alpha_{4}, \alpha_{6}] = \alpha_{3}, [\alpha_{4}, \alpha_{5}] = \alpha_{2}, [\alpha_{3}, \alpha_{6}] = \alpha_{1}, [\alpha_{2}, \alpha_{5}] =  \alpha_{1}^{-1},\\
			& \alpha_{4}^{p} = \alpha_{1}, \alpha_{5}^{p} = \alpha_{3}\alpha_{1}, \alpha_{6}^{p} = \alpha_{2}\alpha_{1}^{-1},  \alpha_{1}^{p} = \alpha_{2}^{p} = \alpha_{3}^{p} =  1 \rangle   \\
			G_{(42,2)}= & \langle \alpha_{1}, \alpha_{2}, \alpha_{3}, \alpha_{4}, \alpha_{5}, \alpha_{6} ~|~ [\alpha_{5}, \alpha_{6}] = \alpha_{4}, [\alpha_{4}, \alpha_{6}] = \alpha_{3}, [\alpha_{4}, \alpha_{5}] = \alpha_{2}, [\alpha_{3}, \alpha_{6}] = \alpha_{1}, [\alpha_{2}, \alpha_{5}] =  \alpha_{1}^{-1},\\
			& \alpha_{4}^{p} = \alpha_{1}, \alpha_{5}^{p} = \alpha_{3},  \alpha_{6}^{p} = \alpha_{2},  \alpha_{1}^{p} = \alpha_{2}^{p} = \alpha_{3}^{p} =  1 \rangle \\
			G_{(42,3k)}= & \langle \alpha_{1}, \alpha_{2}, \alpha_{3}, \alpha_{4}, \alpha_{5}, \alpha_{6} ~|~ [\alpha_{5}, \alpha_{6}] = \alpha_{4}, [\alpha_{4}, \alpha_{6}] = \alpha_{3}, [\alpha_{4}, \alpha_{5}] = \alpha_{2}, [\alpha_{3}, \alpha_{6}] = \alpha_{1}, [\alpha_{2}, \alpha_{5}] =  \alpha_{1}^{-1},\\
			& \alpha_{4}^{p} = \alpha_{1}, \alpha_{5}^{p} = \alpha_{3}\alpha_{1}^{-a-1}, \alpha_{6}^{p} = \alpha_{2}\alpha_{1}^{1-b},  \alpha_{1}^{p} = \alpha_{2}^{p} = \alpha_{3}^{p} =  1 \rangle \text{ where $a$ and $b$ are smallest}\\
			& \text{  positive integers satisfying } a^2 - b^2 \equiv k \bmod p, \text{ for } k=1,2,\ldots,p-1 \\
			G_{(43,1)} = & \langle \alpha_{1}, \alpha_{2}, \alpha_{3}, \alpha_{4}, \alpha_{5}, \alpha_{6} ~|~ [\alpha_{5}, \alpha_{6}] = \alpha_{4}, [\alpha_{4}, \alpha_{6}] = \alpha_{3}, [\alpha_{4}, \alpha_{5}] = \alpha_{2}, [\alpha_{3}, \alpha_{6}] = \alpha_{1}, [\alpha_{2}, \alpha_{5}] =  \alpha_{1}^{-\nu^{-1}}, \\
			& \alpha_{4}^{p} = \alpha_{1}, \alpha_{5}^{p} = \alpha_{3}\alpha_{1}^{-1}, \alpha_{6}^{p} = \alpha_{2}^{\nu}\alpha_{1},  \alpha_{1}^{p} = \alpha_{2}^{p} = \alpha_{3}^{p} =  1 \rangle\\
			G_{(43,2k)} = & \langle \alpha_{1}, \alpha_{2}, \alpha_{3}, \alpha_{4}, \alpha_{5}, \alpha_{6} ~|~ [\alpha_{5}, \alpha_{6}] = \alpha_{4}, [\alpha_{4}, \alpha_{6}] = \alpha_{3}, [\alpha_{4}, \alpha_{5}] = \alpha_{2}, [\alpha_{3}, \alpha_{6}] = \alpha_{1}, [\alpha_{2}, \alpha_{5}] =  \alpha_{1}^{-\nu^{-1}}, \\
			& \alpha_{4}^{p} = \alpha_{1}, \alpha_{5}^{p} = \alpha_{3}\alpha_{1}^{-a-1}, \alpha_{6}^{p} = \alpha_{2}^{\nu}\alpha_{1}^{1-b},  \alpha_{1}^{p} = \alpha_{2}^{p} = \alpha_{3}^{p} =  1 \rangle \text{ where $a$ and $b$ are smallest}\\
			& \text{positive integers satisfying } a^2 - \nu^{-1} b^2 \equiv k \bmod p, \text{ for } k=1,2,\ldots,p-1 \text{ and}\\
			& \nu \text{ be the smallest positive integer which is a non-quadratic residue modulo } p.
		\end{align*}
		\noindent Suppose $G\in \{ G_{(42,1)}, G_{(42,2)}, G_{(42,3k)}, G_{(43,1)}, G_{(43,2k)} \}$. Then $\cd(G) = \{ 1, p, p^2 \}$, $Z(G) = \langle \alpha_{1} \rangle \cong C_{p}$ and $G' = \langle \alpha_{4}, \alpha_{3}, \alpha_{2} \rangle \cong C_{p^2} \times C_{p} \times C_{p}$. 
		For the sake of neatness, we take ${G}' = A$. Then $A$ is a maximal abelian  normal subgroup of $G$. Let $X_G = \{ \chi \}$, where $\chi \in \Irr(G)$ is faithful. Since $G$ is normally monomial, $\chi = \lambda\ind_{A}^{G}$, for some $\lambda \in \lin(A)$. Here, $\ker(\chi) = \Core_{G}(\ker(\lambda)) = 1$. Since $G$ is metabelian, we get $ \left| \frac{A}{\ker(\lambda)} \right| = \exp(A) = p^2$. Thus, $\frac{A}{\ker(\lambda)} \cong C_{p^2}$, and hence, $|\Gamma(\lambda)| = \phi(p^2)$. Through routine computation, it is easy to prove that $\mathbb{Q}(\chi) = \mathbb{Q}(\lambda)$. Thus, $d(\chi) = \chi(1) |\Gamma(\chi)| = p^2 \phi(p^2) = p^3(p-1)$. From Lemma \ref{L1}, we get $m(\chi) = p^3$, and thus, $c(G) = p^4$. \qed

		\subsection{Proof of Theorem \ref{prop:elementaryabelian}} \label{subsec:metabelian}

		Let $G$ be a non-abelian $p$-group of order $p^n$ ($p\geq 3$) such that  $\exp(G) \in \{p, p^2\}$, $\cd(G) = \{ 1, p, p^2 \}$, $d(Z(G)) \geq 2$, and $G$ is not a direct product of an abelian and a non-abelian subgroup. From the hypothesis of Theorem \ref{prop:elementaryabelian}, it is easy to see that $G$ is a metabelian $p$-group.
		Let $d(Z(G)) = m$.	By Lemma \ref{L1}, $|X_{G}|=m$. Suppose $X_G = \{ \chi_{i} \}_{i=1}^{m} \subset \Irr(G)$. Suppose $A$ is a maximal abelian subgroup of $G$ containing $G'$. Then from  \cite[Theorem 2]{BKP}, for $1\leq i \leq m$, there exists $\rho_{i} \in \lin(K_{D_{i}})$, where $D_{i} \leq A$ with $A/D_{i}$ cyclic and $\ker(\rho_{i}\restr_{A}) = D_{i}$, such that 
		\begin{equation} \label{eq:elemabelian}
			\chi_{i} = \rho_{i}\ind_{K_{D_{i}}}^{G} \in \Irr(G),
		\end{equation}
	where for each $i$, $K_{D_i}$ is a 
	fixed maximal element of $\{ T ~|~ A\leq T \leq G \text{ and } T' \leq D_i \}$.
		Since $X_G$ satisfies \eqref{eq:X_G}, we get $\bigcap_{i=1}^{m} \ker\left(\rho_{i}\ind_{K_{D_{i}}}^{G}\right)  = 1$.  From \cite[Theorem 12.11]{I}, $G$ does not have any abelian subgroup of index $p$ in $G$. In this case, we take $A$ to be an elementary abelian subgroup of index $p^2$ in $G$. Suppose $i=1$. We have three possibilities for $K_{D_{1}}$, namely, $G$, $K$, and $A$, where $K \leq G$ such that $|G/K|=p$ and $A\subset K$ with $K{}' \subseteq D_1$. \\
		\noindent	{\bf Case I ($K_{D_{1}} = G$):} Here $\frac{G}{\ker(\rho_{1})} \cong C_{p^2}$ or $C_{p}$. Now, if $\frac{G}{\ker(\rho_{1})} \cong C_{p}$, then $d(\rho_{1}) = \phi(p)$. Then, from Lemma \ref{L1}, we get that $c(G) = \delta_1 + \delta_2 p$, for some $\delta_1, \delta_2 \in \mathbb{N}$, where $\delta_2 \neq 0$. Since $G$ is not a direct product of an abelian and a non-abelian subgroup, we get a contradiction from \cite[Theorem 3.2]{BG} and Lemma \ref{lemma:a_1=0}. 
		Hence, $\frac{G}{\ker(\rho_{1})} \cong C_{p^2}$, and so $d(\rho_{1}) = \phi(p^2) = p(p-1)$.\\
		{\bf Case II ($K_{D_{1}} = K$):} Here $\frac{K}{\ker(\rho_{1})} \cong C_{p^2}$ or $C_{p}$. Then, $d(\rho_{1}\ind_{K}^{G}) = p |\Gamma(\rho_{1}\ind_{K}^{G})| \leq p |\Gamma(\rho_{1})|\leq p\cdot p(p-1) = p^2(p-1)$. On the other hand, $d(\rho_{1}\ind_{K}^{G}) \geq \rho_{1}\ind_{K}^{G}(1)\cdot (p-1)= p(p-1)$. Therefore, $p(p-1) \leq d(\rho_{1}\ind_{K}^{G}) \leq p^2(p-1)$. Then we get $d(\rho_{1}\ind_{K}^{G}) = p(p-1)$, or $p^2(p-1)$.\\
		{\bf Case III ($K_{D_{1}} = A$):} In this case, we get $\frac{A}{\ker(\rho_{1})} \cong C_{p}$, and $d(\rho_{1}) = p-1$. This implies that $d(\rho_{1}\ind_{A}^{G}) \leq p^2(p-1)$. On the other hand, $d(\rho_{1}\ind_{A}^{G}) \geq p^2(p-1)$. Hence, we get $d(\rho_{1}\ind_{A}^{G}) = p^2(p-1)$.\\
		Similarly, we get $d(\rho_{i}\ind_{K_{D_{i}}}^{G}) = p(p-1)$, or $p^2(p-1)$, for $2\leq i\leq m$. Thus, $\sum_{i=1}^{m}d(\rho_{i}\ind_{K_{D_{i}}}^{G}) = ap(p-1) + bp^2(p-1)$, for some $0 \leq a, b < p$ such that $a+b = m$. Then from Lemma \ref{L1}, we get $c(G) = ap^2 + bp^3$, for some $0 \leq a, b < p$ such that $a+b = m$. \qed \\
		
		\noindent Here we present examples of some $p$-groups $G$ such that $d(Z(G))=2$ and $c(G)$ takes all the possible values obtained in Theorem \ref{prop:elementaryabelian}. 
		
		\begin{example} \label{example:elementaryabelian} \textnormal{In the following group presentations, all relations of the form $[\alpha, \beta] = 1$ (with $\alpha$, $\beta$ generators) have been omitted.} 	
			\begin{enumerate}
				\item  \textnormal{Consider the groups
					\begin{align*}
						G_{1} &= \phi_{4}(221)b = \langle \alpha, \alpha_{1}, \alpha_{2}, \beta_{1}, \beta_{2} ~|~ [\alpha_{i}, \alpha] = \beta_{i}, \alpha^{p} = \beta_{2}, \alpha_{2}^{p} = \beta_{1}, \alpha_{1}^{p} = \beta_{i}^{p} = 1 ~ (i=1,2) \rangle, \\
						G_{2} &= \phi_{4}(221)f_{0} = \langle \alpha, \alpha_{1}, \alpha_{2}, \beta_{1}, \beta_{2} ~|~ [\alpha_{i}, \alpha] = \beta_{i}, \alpha_{1}^{p} = \beta_{2}, \alpha_{2}^{p} = \beta_{1}^{\nu}, \alpha^{p} = \beta_{i}^{p} = 1 ~ (i=1,2) \rangle, \\
						\text{and}& ~ G_{3} = \phi_{4}(2111)a = \langle \alpha, \alpha_{1}, \alpha_{2}, \beta_{1}, \beta_{2} ~|~ [\alpha_{i}, \alpha] = \beta_{i}, \alpha^{p} = \beta_{2}, \alpha_{i}^{p} = \beta_{i}^{p} = 1 ~ (i=1,2) \rangle
					\end{align*}
					of order $p^5$ ($p\geq 5$) belonging to the isoclinic family $\Phi_{4}$ (see \cite[Section 4.5]{RJ}). Then for $1\leq i \leq 3$, $\exp(G_{i}) = p^2$, $Z(G_i) = \langle \beta_{1}, \beta_{2} \rangle$, $d(Z(G_i)) = 2$, $\cd(G_{i}) = \{ 1, p \}$ and $G_{i}$ is not a direct product of an abelian and a non-abelian subgroup. From the presentation of $G_{i}$, for $1\leq i \leq 3$, it is easy to see that $G_{i}$ has an abelian normal subgroup of index $p$, namely, $K=\langle \alpha_{1}, \alpha_{2}, \beta_{1}, \beta_{2} \rangle $. From Remark \ref{remark:nMcd(G)1,p}, we have $c(G_i) = ap^2 + bp^3$, for some $0 \leq a, b < p$ such that $a+b = 2$. Then $c(G_{i}) \in \{ 2p^2, p^3+p^2, 2p^3 \}$, for each $i$. It is easy to check that $c(G_{1}) = p^3+p^2$, $c(G_{2}) = 2p^3$ and $c(G_{3}) = 2p^2$.	}
				\item \textnormal{ Consider the groups
					\begin{align*}
						G_{1}  = \langle \alpha_{1}, \alpha_{2}, \alpha_{3}, \alpha_{4}, \alpha_{5}, \alpha_{6} ~|~& [\alpha_{3}, \alpha_{4}] = \alpha_{1}, [\alpha_{5}, \alpha_{6}] = \alpha_{2}, \alpha_{3}^{p} = \alpha_{2},\\
						& \alpha_{1}^{p} = \alpha_{2}^{p} =  \alpha_{4}^{p} = \alpha_{5}^{p}= \alpha_{6}^{p} =  1 \rangle, \\
						G_{2}  = \langle \alpha_{1}, \alpha_{2}, \alpha_{3}, \alpha_{4}, \alpha_{5}, \alpha_{6} ~|~& [\alpha_{3}, \alpha_{4}] = \alpha_{1}, [\alpha_{5}, \alpha_{6}] = \alpha_{2}, \alpha_{3}^{p} = \alpha_{1},\\
						& \alpha_{1}^{p} = \alpha_{2}^{p} =  \alpha_{4}^{p} = \alpha_{5}^{p}= \alpha_{6}^{p} =  1 \rangle, \text{ and}\\
						G_{3} =  \langle \alpha_{1}, \alpha_{2}, \alpha_{3}, \alpha_{4}, \alpha_{5}, \alpha_{6} ~|~& [\alpha_{3}, \alpha_{4}] = \alpha_{1}, [\alpha_{5}, \alpha_{6}] = \alpha_{2}, \alpha_{4}^{p} = \alpha_{5}^{p} = \alpha_{1}\alpha_{2},\\
						& \alpha_{1}^{p} = \alpha_{2}^{p} =  \alpha_{3}^{p} =  \alpha_{6}^{p} =  1 \rangle
					\end{align*}
					of order $p^6$ ($p\geq 5$) belonging to the isoclinic family $\Phi_{12}$ (see \cite{Easterfield}). Then for $1\leq i \leq 3$, $|G_{i}| = p^6$, $\exp(G_{i}) = p^2$, $Z(G_i) = \langle \alpha_{1}, \alpha_{2} \rangle$, $d(Z(G_i)) = 2$, $\cd(G_{i}) = \{ 1, p, p^2 \}$ and $G_{i}$ is not a direct product of an abelian and a non-abelian subgroup. It is easy to see that $\langle \alpha_{1},\alpha_{2}, \alpha_{4}, \alpha_{6} \rangle$ and  $\langle \alpha_{1}, \alpha_{2}, \alpha_{3}, \alpha_{6}\rangle$ are elementary abelian normal subgroups of index $p^2$ in $G_i$ ($i=1,2$) and $G_3$, respectively. From Theorem \ref{prop:elementaryabelian}, $c(G_i) = ap^2 + bp^3$, for some $0 \leq a, b < p$ such that $a+b = 2$. Then $c(G_{i}) \in \{ 2p^2, p^3+p^2, 2p^3 \}$, for each $i$. It is easy to check that $c(G_{1}) = p^3+p^2$, $c(G_{2}) = 2p^2$ and $c(G_{3}) = 2p^3$.	}
			\end{enumerate}
		\end{example}
		
		\begin{remark} \label{remark:elementaryabelianexample}
			\textnormal{Let $G$ be a non-abelian $p$-group of order $p^n$ ($p\geq 3$) such that  $\exp(G) \in \{ p, p^2 \}$, $d(Z(G)) \geq 2$, $\cd(G) = \{ 1, p, p^2 \}$ and $G$ is not a direct product of an abelian and a non-abelian subgroup. Then $G$ may not have an elementary abelian subgroup of index $p^2$ in $G$. For example, for $p\geq 5$, consider
				\begin{align*}
					G =	\langle \alpha_{1}, \ldots, \alpha_{6} ~|~ [\alpha_{5}, \alpha_{6}] = \alpha_{3}, & [\alpha_{4}, \alpha_{5}] = \alpha_{2},  [\alpha_{3}, \alpha_{6}] =  \alpha_{1}, \alpha_{4}^{p} = \alpha_{1}, \alpha_{5}^{p}= \alpha_{2}, \alpha_{1}^{p} = \alpha_{2}^{p} = \alpha_{3}^{p} =  \alpha_{6}^{p} = 1 \rangle,
				\end{align*}
				which is a $p$-group of order $p^6$ belonging to the isoclinic family $\Phi_{17}$ (see \cite{Easterfield}). Suppose $A$ is a normal abelian subgroup of index $p^2$ in $G$. Then $G'$ is contained in $A$. Through routine computation, it is easy to see that $\exp(A)$ must be $p^2$.
			}
		\end{remark}

		\section{Acknowledgements}
		Ayush Udeep acknowledges University Grants Commission, Government of India (File no: Nov2017-434175). The corresponding author acknowledges SERB, Government of India for financial support through grant (MTR/2019/000118). The authors are grateful to the anonymous referee for his/her valuable suggestions and for providing a shorter version of some of the proofs. He also provided the corrected proof of \cite[Lemma 2.2]{BG}.

	\end{document}